\documentclass{article}
\usepackage{graphicx}
\def\be{\begin{equation}}
\def\ee{\end{equation}}
\newcommand{\ff}[1]{{\mbox{\boldmath $#1$}}}
\def\a{\alpha}
\def\b{\beta}
\def\x{\ff{x}}
 
\begin{document}

\title{\large \bf Multiobjective Firefly Algorithm for Continuous Optimization}

\author{Xin-She Yang \\
Department of Engineering,  University of Cambridge, \\
Trumpington Street, Cambridge CB2 1PZ, UK }

\date{}

\maketitle

\begin{abstract}
Design problems in industrial engineering often involve a large number of
design variables with multiple objectives, under complex nonlinear constraints.
The algorithms for multiobjective problems can be significantly
different from the methods for single objective optimization.
To find the Pareto front and non-dominated set for a nonlinear
multiobjective optimization problem may require significant computing effort,
even for seemingly simple problems.
Metaheuristic algorithms start to show their advantages in
dealing with multiobjective optimization. In this paper, we extend the
recently developed firefly algorithm to solve multiobjective optimization
problems. We validate the proposed approach using a selected subset of
test functions and then apply it to solve design optimization benchmarks.
We will discuss our results and provide topics for further research.  \\


{\it Citation details:} X. S. Yang, Multiobjective firefly algorithm for continuous optimization,
{\it Engineering with Computers}, Vol. {\bf 29}, Issue 2, pp. 175--184 (2013).

\end{abstract}

\bigskip

{\it Revised manuscript: Ms. No. EWCO-D-11-00145 (Engineering with Computers) }

\newpage


\section{Introduction}

Design optimization in engineering and industry often concerns multiple design objectives
under complex, highly nonlinear constraints. Different objectives often conflict
each other, and sometimes, truly optimal solutions do not exist, and some
compromises and approximations are often needed \cite{Cag,Deb2,Leif}. Further to this complexity,
a design problem is subjected to various design constraints, limited by
design codes or standards, material properties and the optimal utility of available resources
and costs \cite{Deb2,Farina}. Even for global optimization problems with a single objective, if the design
functions are highly nonlinear, global optimality is not easy to reach.
Metaheuristic algorithms are very powerful in dealing with this kind of
optimization, and there are many review articles and textbooks
\cite{Coello,Deb,Geem,Talbi,Yang,Yang2,Yang3}.

In contrast with single objective optimization, multiobjective problems are much more
difficult and complex \cite{Coello,Gong}. Firstly, no single unique solution is the best;
instead, a set of non-dominated solutions should be found in order to get a good approximation
to the true Pareto front. Secondly, even if an algorithm can find solution points
on the Pareto front, there is no guarantee that multiple Pareto points will distribute
along the front uniformly, often they do not. Thirdly, algorithms which work well
for single objective optimization usually cannot directly work for multiobjective problems,
unless under special circumstances such as combining multiobjectives into a single
objective using some weighted sum method. Substantial modifications are needed
to make algorithms for single objective optimization work.
In addition to these difficulties, a further challenge
is how to generate solutions with enough diversity so that new solutions can
sample the search space efficiently.

Furthermore, real-world optimization problems always involve some degree of
uncertainty or noise. For example, materials properties for a design product
may vary significantly. An optimal design should be robust enough to allow
such inhomogeneity, which provides a set of multiple feasible solution sets.
Consequently, optimal solutions among the robust Pareto set can provide
good options so that decision-makers or designers can choose to suit their needs.
Despite these challenges, multiobjective optimization has many powerful algorithms
with many successful applications \cite{Deb,Abbass,Banks,Konak,Rang}. In addition, metaheuristic algorithms start to emerge as
a major player for multiobjective global optimization, they often mimic the
successful characteristics in Nature, especially biological systems \cite{Yang,Yang2},
while some algorithms are inspired by the beauty of music \cite{Geem01}.
Many new algorithms are emerging with many important
applications \cite{Talbi,Yang2,Abbass,Ken,Osy,Reyes,ZhangQ}.

For example, multiobjective genetic algorithms are widely known \cite{Konak,Schaf}, while multiobjective differential evolution
algorithms are also very powerful \cite{Robic,Xue}. In addition, multiobjective
particle swarm optimizers are becoming increasingly popular \cite{Reyes}.
As there are many algorithms, one of our motivations in the present study is to compare
the performance of these algorithms for real-world application.

Most metaheuristic algorithms are based on the so-called swarm intelligence.
PSO is a good example, it mimics some characteristics of birds and fish swarms.
Recently, a new metaheuristic search algorithm, called Firefly Algorithm (FA),
has been developed by Yang \cite{Yang,Yang3}.  FA mimics some characteristics of
tropic firefly swarms and their flashing behaviour \cite{Yang2,Yang3}. A firefly tends
to be attracted towards other fireflies with higher flash intensity. This algorithm is thus
different from PSO and can have two advantages: local attractions and automatic regrouping.
As light intensity decreases with distance, the attraction among fireflies can be local
or global, depending on the absorbing coefficient, and thus all local modes as well as
global modes will be visited. In addition, fireflies can also subdivide and thus regroup into
a few subgroups due to neighboring attraction is stronger than long-distance attraction,
thus it can be expected each subgroup will swarm around a local mode. This latter advantage
makes it particularly suitable for multimodal global optimization problems \cite{Yang2,Yang3}.

Preliminary studies show that
it is very promising and could outperform existing algorithms such as particle swarm optimization (PSO).
For example, a Firefly-LGB algorithm, based on firefly algorithm and Linde-Buzo-Gray (LGB)
algorithms for vector quantization of digital image compression, was developed by
Horng and Jiang \cite{Horng}, and their results suggested that Firefly-LGB is faster
than other algorithms such as particle swarm optimization LBG (PSO-LBG) and
honey-bee mating optimization LBG (HBMO-LBG).
Apostolopoulos and Vlachos provided a detailed background and analysis over a wide range of
test problems \cite{Apost}, and they also solved multiobjective load dispatch problem using a weighted sum method
by combining multiobjectives into a single objective, and their results are very promising.
The preliminary successful results of firefly-based algorithms provide another motivation
for this paper. That is to see how this algorithm can be extended to solve multiobjective
optimization problems.

In this paper, we will extend FA to solve multiobjective problems and formulate
a multiobjective firefly algorithm (MOFA). We will first validate it
against a subset of multiobjective test functions.
Then, we will apply it to solve design optimisation problems
in engineering, including bi-objective beam design and a design of a disc brake.
Finally, we will discuss the unique features of the proposed algorithm as well as topics for further studies.

\section{Multiobjective Firefly Algorithm}

In order to extend the firefly algorithm for single objective
optimization to solve multiobjective problems, let us briefly review its
basic version.

\subsection{The Basic Firefly Algorithm}

Firefly Algorithm  was developed by Yang for continuous optimization \cite{Yang,Yang2,YangFA2},
which was subsequently applied into structural optimization \cite{Gandomi} and image processing \cite{Horng}.
FA was based on the flashing patterns and behaviour
of fireflies. In essence, FA uses the following three idealized rules: (1)
Fireflies are unisex so that one firefly will be attracted to other fireflies
regardless of their sex; (2) The attractiveness of a firefly is proportional to
its brightness and they both decrease with distance. Thus for any two flashing fireflies,
the less brighter one will move towards the brighter one.  If there is
no brighter one than a particular firefly, it will move randomly;
(3) The brightness of a firefly is determined by the landscape of the
objective function.

For a  maximization problem, the brightness can simply be proportional
to the value of the objective function.
As both light intensity and attractiveness affect the movement of fireflies
in the firefly algorithm, we have to define their variations. For simplicity,
we can always assume that the attractiveness of a firefly is
determined by its brightness which in turn is associated with the
encoded objective function. In the simplest case for maximum optimization problems,
the brightness $I$ of a firefly at a particular location $\x$ can be chosen
as $I(\x) \propto f(\x)$.

However, the attractiveness $\b$ is relative, it should be
seen in the eyes of the beholder or judged by the other fireflies. Thus, it will
vary with the distance $r_{ij}$ between firefly $i$ and firefly $j$. Therefore,
we can now define the attractiveness $\b$ of a firefly by
\be \b = \b_0 e^{-\gamma r^2}, \label{att-equ-100} \ee
where $\b_0$ is the attractiveness at $r=0$.
In fact, equation (\ref{att-equ-100}) defines a characteristic distance
$\Gamma=1/\sqrt{\gamma}$ over which the attractiveness changes significantly
from $\b_0$ to $\b_0 e^{-1}$.
The distance between any two fireflies $i$ and $j$ at $\x_i$ and $\x_j$, respectively, is
the Cartesian distance $ r_{ij}=||\x_i-\x_j||$.
It is worth pointing out that the distance $r$ defined above is {\it not} limited
to the Euclidean distance. In fact, any measure that can effectively characterize
the quantities of interest in the optimization problem can be used as the `distance' $r$.
We can define other distance $r$ in the $n$-dimensional
hyperspace, depending on the type of problem of our interest.

For any given two fireflies $\x_i$ and $\x_j$,
the movement of firefly $i$ is attracted to another more attractive (brighter)
firefly $j$ is determined by
\be \x_i^{t+1} =\x_i^t + \b_0 e^{-\gamma r^2_{ij}} (\x_j^t-\x_i^t) + \alpha_t \; \ff{\epsilon}_i^t,
\label{FA-equ-50}  \ee
where the second term is due to the attraction. The third term
is randomization with $\alpha_t$ being the randomization parameter, and
$\ff{\epsilon}_i^t$ is a vector of random numbers drawn from a Gaussian distribution
or uniform distribution. The location of fireflies can be updated sequentially, by comparing and
updating each pair of them in every iteration cycle.

For most implementations,
we can take $\b_0=1$ and $\alpha_t=O(1)$, though we found that it is better to use
a time-dependent $\alpha_t$ so that randomness can be reduced gradually as iterations proceed.
It is worth pointing out that (\ref{FA-equ-50}) is a random walk biased
towards the brighter fireflies. If $\b_0=0$, it becomes a simple random walk.
Furthermore, the randomization term can easily be
extended to other distributions such as L\'evy flights \cite{YangDeb}.

The parameter $\gamma$ now characterizes the variation of the attractiveness,
and its value is crucially important in determining the speed of the convergence
and how the FA algorithm behaves. In theory, $\gamma \in [0,\infty)$, but in
practice, $\gamma=O(1)$ is determined by the characteristic distance $\Gamma$ of the
system to be optimized. Thus, for most applications, it typically varies
from $10^{-5}$ to $10^{5}$.

To consider the scale variations of each problem, we now use rescaled, vectorized parameters
\be \a=0.01 L, \quad \gamma=0.5/L^2, \ee
with $L=(U_b-L_b)$ where $U_b$ and $L_b$ are the upper and lower bounds of $\x$, respectively.
Here the factor $0.01$ is to make
sure the random walks is not too aggressive, and this value has been obtained by
a parametric study.

\subsection{Multiobjective Firefly Algorithm}

For multiobjective optimization, one way is to combine all objectives into a
single objective so that algorithms for single objective optimization can
be used without much modifications. For example, FA can be used directly to solve
multiobjective problems in this manner, and a detailed study was carried out
by Apostolopoulos and Vlachos \cite{Apost}.

Another way is to extend the firefly algorithm to produce
Pareto optimal front directly. By extending the basic ideas of FA, we can
develop the following Multi-objective Firefly Algorithm (MOFA), which
can be summarized as the pseudo code listed in Fig. \ref{code-40}.

\begin{figure}[h]
\begin{center}
\begin{minipage}[c]{0.9\textwidth}
\hrule \vspace{5pt}
\indent \quad Define objective functions $f_1(\x),..., f_K(\x)$ where $\x=(x_1, ..., x_d)^T$ \\
\indent \quad Initialize a population of $n$ fireflies $\x_i \; (i=1,2,...,n)$ \\
\indent \quad {\bf while} ($t<$MaxGeneration) \\
\indent \qquad {\bf for} $i,j=1:n$ (all $n$ fireflies) \\
\indent \qquad \quad Evaluate their approximations $PF_i$ and $PF_j$ to the Pareto front\\
\indent \qquad \qquad \qquad if $i \ne j$ and when all the constraints are satisfied \\
\indent \qquad \quad {\bf if} $PF_j$ dominates $PF_i$, \\
\indent \qquad \qquad  Move firefly $i$ towards $j$ using (\ref{FA-equ-50}) \\
\indent \qquad \qquad Generate new ones if the moves do not satisfy all the constraints \\
\indent \qquad \quad {\bf end if} \\
\indent \qquad \quad  {\bf if} no non-dominated solutions can be found \\
\indent \qquad \qquad Generate random weights $w_k$ ($k=1,...,K$) \\
\indent \qquad \qquad Find the best solution $\ff{g}_*^t$ (among all fireflies) to minimize $\psi$ in (\ref{mofa-equ-100}) \\
\indent \qquad \qquad Random walk around $\ff{g}_*^t$  using (\ref{mofa-equ-150}) \\
\indent \qquad \quad {\bf end if} \\
\indent \qquad  \quad Update and pass the non-dominated solutions to next iterations \\
\indent \qquad {\bf end} \\
\indent \qquad Sort and find the current best approximation to the Pareto front  \\
\indent \qquad Update  $t \leftarrow t+1$ \\
\indent \quad {\bf end while} \\
\indent \quad Postprocess results and visualisation;
\hrule \vspace{5pt}
\caption{Pseudo Code: Multiobjective firefly algorithm (MOFA). \label{code-40} }
\end{minipage}
\end{center}
\end{figure}

The procedure starts with an appropriate definition of objective functions with associated
nonlinear constraints. We first initialize a population of $n$ fireflies so that they should distribute
among the search space as uniformly as possible. This can be achieved by using sampling techniques
via uniform distributions. Once the tolerance or a fixed number of iterations is defined,
the iterations start with the evaluation of brightness or objective values of all the fireflies
and compare each pair of fireflies. Then, a random weight vector is generated (with the sum equal to 1),
so that a combined best solution $\ff{g}_*^t$ can be obtained. The non-dominated solutions are
then passed onto the next iteration. At the end of a fixed number of iterations, in general $n$ non-dominated
solution points can be obtained to approximate the true Pareto front.

In order to do random walks more efficiently, we can find the current best
$\ff{g}_*^t$ which minimizes a combined objective via the weighted sum
\be \psi(\x) =\sum_{k=1}^K w_k f_k, \quad \sum_{k=1}^K w_k=1. \label{mofa-equ-100} \ee
Here $w_k=p_k/K$ where $p_k$ are the random numbers drawn from a uniform distributed Unif[0,1].
In order to ensure that $\sum_k w_k=1$, a rescaling operation is performed after generating
$K$ uniformly distributed numbers.
It is worth pointing out that the weights $w_k$ should be chosen randomly at each
iteration, so that the non-dominated solution can sample diversely along the
Pareto front.

If a firefly is not dominated by others in the sense of Pareto front,
the firefly moves
\be \x_i^{t+1} =\ff{g}_*^t + \alpha_t  \; \ff{\epsilon}_i^t,
\label{mofa-equ-150}  \ee
where $\ff{g}_*^t$ is the best solution found so far for a given set of random weights.

Furthermore, the randomness can be reduced as the iterations proceed, and this can be
achieved in a similar manner as that for simulated annealing and other random reduction
techniques \cite{Yang3}. We will use
\be \a_t=\a_0 0.9^t, \ee
where $\a_0$ is the initial randomness factor.

\subsection{Pareto Optimal Front}

For a minimization problem, a solution vector $\ff{u}=(u_1, .., u_n)^T$
is said to dominate another vector
$\ff{v}=(v_1,..., v_n)^T$ if and only if $u_i \le v_i$ for $\forall i \in \{ 1,...,n\}$
and $\exists i \in \{1,...,n \}: u_i < v_i. $
In other words, no component of $\ff{u}$ is larger than the corresponding component of $\ff{v}$,
and at least one component is smaller.
Similarly, we can define another dominance relationship $\preceq$ by
\be \ff{u} \preceq \ff{v} \Longleftrightarrow \ff{u} \prec \ff{v} \vee \ff{u}=\ff{v}. \ee
It is worth pointing out that for maximization problems, the dominance can be defined
by replacing $\prec$ with $\succ$. Therefore,
a point $\x_*$ is called a non-dominated solution if
no solution can be found that dominates it \cite{Coello}.

The Pareto front $PF$ of a multiobjective can be defined as the
set of non-dominated solutions so that
\be  PF =\{ \ff{s} \in {\cal S} \Big|
\exists \hspace{-5pt} / \; \ff{s'} \in {\cal S}: \ff{s'} \prec \ff{s} \}, \ee
where ${\cal S}$ is the solution set.

To obtain a good approximation of the Pareto front, a diverse range of solutions should be
generated using efficient techniques \cite{Konak,Erfani,Guj,Marler}.
For example, L\'evy flights ensure
the good diversity of the solutions, as we can see from later simulations.

\section{Numerical Results}

We have implemented the proposed MOFA in Matlab, and we have first validated it against
a set of multiobjective test functions. Then, we have used it to solve some industrial
design of structures. In order to obtain the right algorithm-dependent parameters,
we have carried out detailed parametric studies.

\subsection{Parametric Studies}

By varying the parameters $\a_0$, $\b_0$ and $\gamma$, we have carried out
parametric studies by setting $\a_0=0$ to $1$ with a step of $0.05$,
$\b_0=0$ to $1$ with a step of $0.05$, and $\gamma=0.1$ to $10$ with a step of $0.1$
and then $1$. From simulations, we found that
we can use $\a_0=0.1$ to $0.5$, $\b_0=0.7$ to $1.0$ and $\gamma=1$
for most problems.

The stopping criterion can be defined in many ways. We can either use
a given tolerance or a fixed number of iterations. For a given tolerance,
we should have some prior knowledge of the true optimum of the objective function
so that we can calculate the differences between the current best solutions
and the true optimal solution so that we can assess if the tolerance is met.
In reality, we usually do not know the true optimum in advance, except for
a few well-tested cases. In addition, the number of
functions can vary significantly from function to function even for the same
tolerance.  From the implementation
point of view, a fixed number of iterations is not only easy to implement,
but also suitable to compare the closeness of Pareto front of different functions.
So we have set the fixed number iterations as 2500, which is sufficient
for most problems. If necessary, we can also increase it to a larger number.

One can generate  points of the Pareto front in two ways: increase the population size $n$ or run the program
a few more times. Through simulations, we found that  increasing $n$ typically
leads to a longer computing time than re-running the program a few times.
This may be due to the fact that manipulations of large matrices or longer vectors
usually take longer. In order to generate $M$ points using a smaller population size $n$,
it requires to run the program $M/n$ times, each run with different, random initial configurations
but with the same number of iterations $t=2500$. For example, to generate $M=200$ points,
we can use $n=50$, which is easily done within a few minutes.
Therefore, in all our simulations, we will use the fixed parameters: $n=50$,
$\a_0=0.25$, $\b_0=1$ and  $\gamma=1$.

\subsection{Multiobjective Test Functions}

\begin{figure}
\centerline{\includegraphics[height=3in,width=4in]{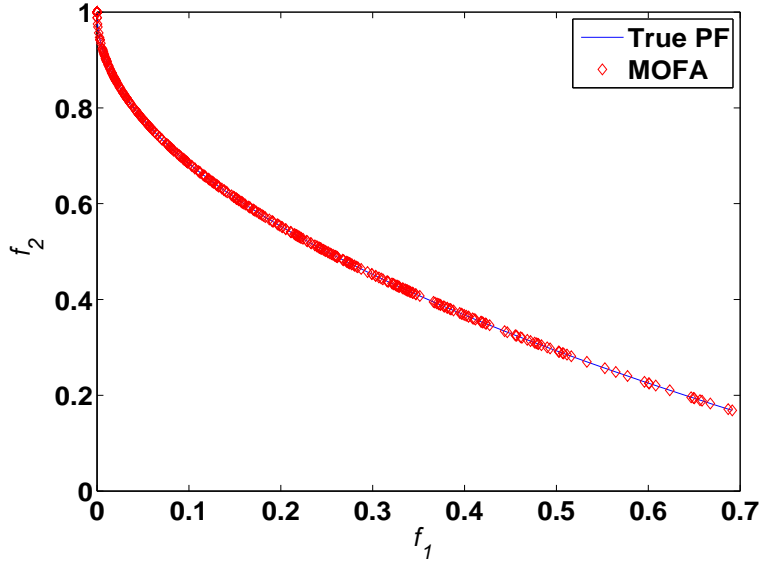} }
\caption{Pareto front of ZDT1:
a comparison of the front found by MOFA and the true Pareto front (true PF).
Here the horizontal axis is $f_1$ while the vertical axis is $f_2$.
\label{fig-100}}
\end{figure}

There are many different test functions for multiobjective optimization \cite{ZhangZhou,Zit,Zitz},
but a subset of a few widely used
functions provides a wide range of diverse properties in terms of Pareto front
and Pareto optimal set. To validate the proposed MOFA, we have selected a subset
of these functions with convex, non-convex and discontinuous Pareto fronts.
We also include functions with more complex Pareto sets. To be more specific
in this paper, we have tested the following 5 functions:

\begin{itemize}
\item Schaffer's Min-Min (SCH) test function with convex Pareto front \cite{Schaf,Zhang}
\be f_1(x)=x^2, \quad f_2(x) =(x-2)^2, \quad -10^{3} \le x \le 10^3. \ee

\item ZDT1 function with a convex front \cite{Zit,Zitz}
\[ f_1(x)=x_1, \quad f_2(x)=g (1-\sqrt{f_1/g}), \]
\be g=1+\frac{9 \sum_{i=2}^d x_i}{d-1}, \quad x_i \in [0,1], \; i=1,...,30, \ee
where $d$ is the number of dimensions. The Pareto optimality is reached when $g=1$,
and thus the true Pareto front is $f_2=1-\sqrt{f_1}$.

\item ZDT2 function with a non-convex front
\[ f_1(x)=x_1, \quad f_2(x) =g (1-\frac{f_1}{g})^2, \]

\item ZDT3 function with a discontinuous front
\[ f_1(x) =x_1, \quad f_2(x)=g \Big[1-\sqrt{\frac{f_1}{g}}-\frac{f_1}{g} \sin (10 \pi f_1) \Big], \]
where $g$ in functions ZDT2 and ZDT3 is the same as in function ZDT1. In the ZDT3 function,
$f_1$ varies from $0$ to $0.852$ and $f_2$ from $-0.773$ to $1$.

\item LZ function \cite{ZhangQ,Li}
\[ f_1=x_1 +\frac{2}{|J_1|} \sum_{j \in J_1} \Big [ x_j -\sin (6 \pi x_1 +\frac{j \pi}{d}) \Big]^2, \]
\be f_2=1-\sqrt{x_1} + +\frac{2}{|J_2|} \sum_{j \in J_2} \Big [ x_j -\sin (6 \pi x_1 +\frac{j \pi}{d}) \Big]^2, \ee
where $J_1=\{j|j$ is odd $\}$ and $J_2 =\{ j|j$ is even $\}$ where $2 \le j \le d$. This
function has a Pareto front $f_2=1-\sqrt{f_1}$ with a Pareto set
\be x_j=\sin (6 \pi x_1 + \frac{j \pi}{d}), \quad j=2,3, ..., d, \quad x_1 \in [0,1]. \ee

\end{itemize}
After generating 200 Pareto points by MOFA, these points are
compared with the true front $f_2=1-\sqrt{f_1}$ of ZDT1 (see Fig. \ref{fig-100}).

Let us define the distance or error between the estimated Pareto front $PF^e$ to
its corresponding true front $PF^t$ as
\be E_f=||PF^e-PF^t||^2=\sum_{j=1}^N (PF_j^e-PF_j^t)^2, \ee
where $N$ is the number of points. For all the test functions, the true Pareto fronts
have analytical forms \cite{ZhangQ,Zit,Zhang}, for example, $f_2=1-\sqrt{f_1}$ for ZDT1,
which makes the above calculations straightforward.

The convergence property can be viewed by
following the iterations. As this measure is an absolute measure, which depends on
the number of points. Sometimes, it is easier to use relative measure using generalized
distance
\be D_g=\frac{1}{N} \sqrt{\sum_{j=1}^N (PF_j-PF_j^t)^2}. \ee

\begin{figure}
\centerline{\includegraphics[height=3in,width=4in]{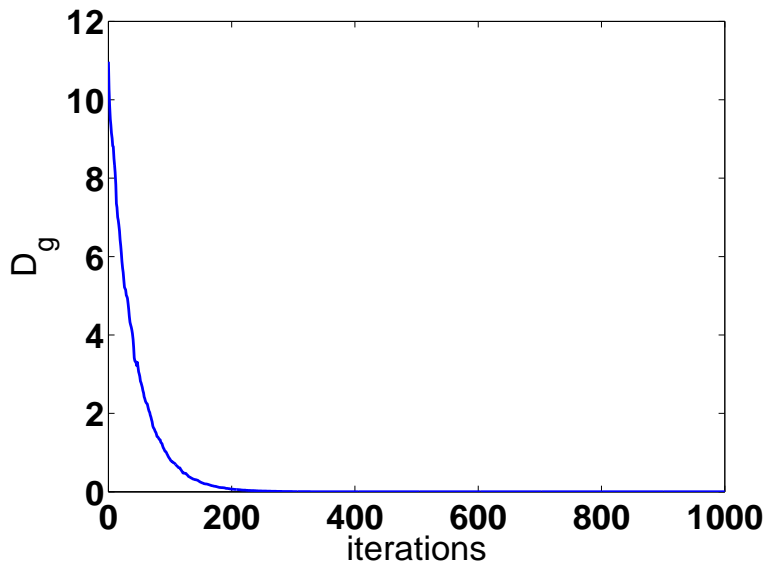} }
\caption{Convergence of the proposed MOFA. The least-square distance (vertical axis)
from the estimated front to the true front of ZDT1 for the first 1000 iterations. \label{fig-200}}
\end{figure}

\begin{figure}
\centerline{\includegraphics[height=2.5in,width=3in]{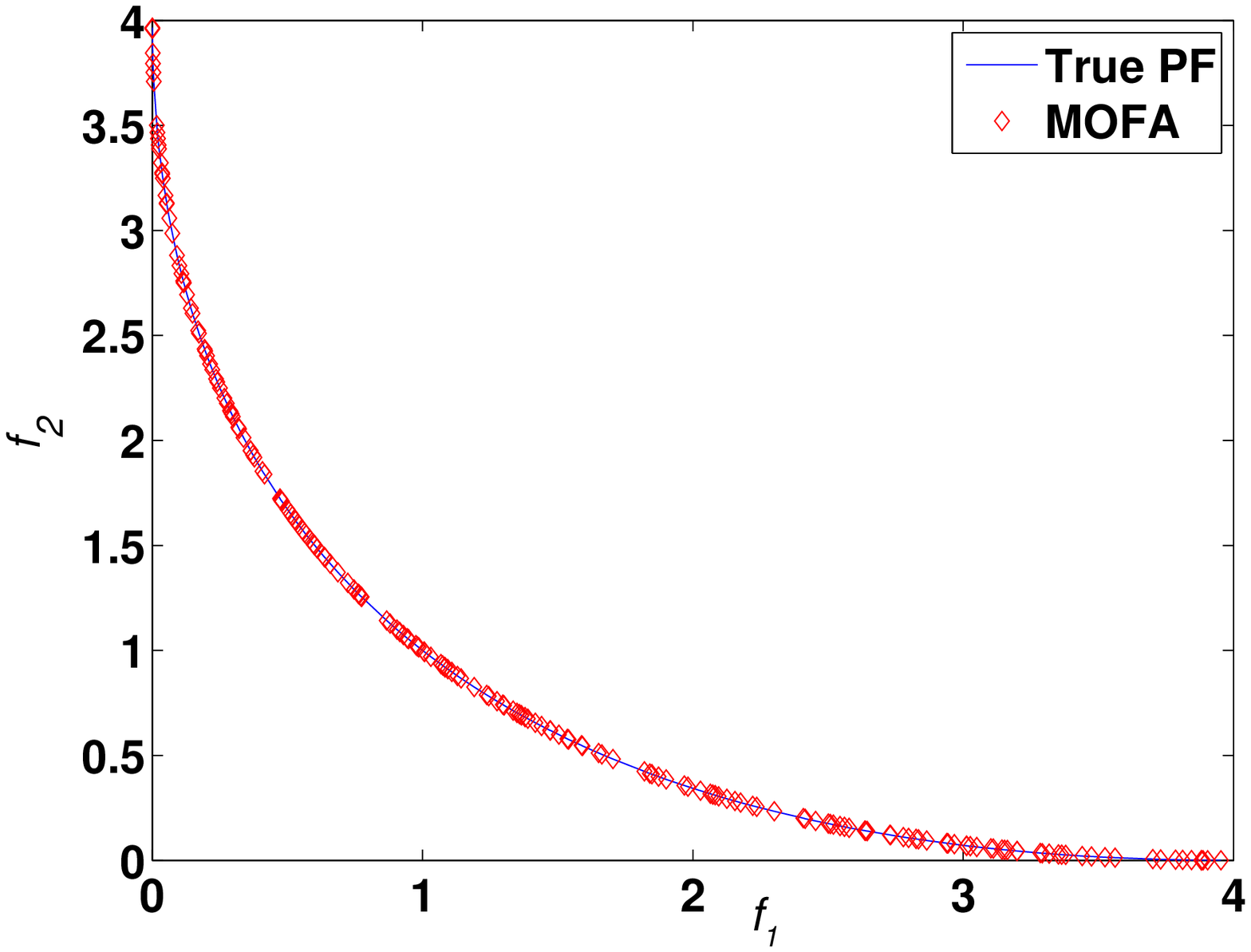}
\includegraphics[height=2.5in,width=3in]{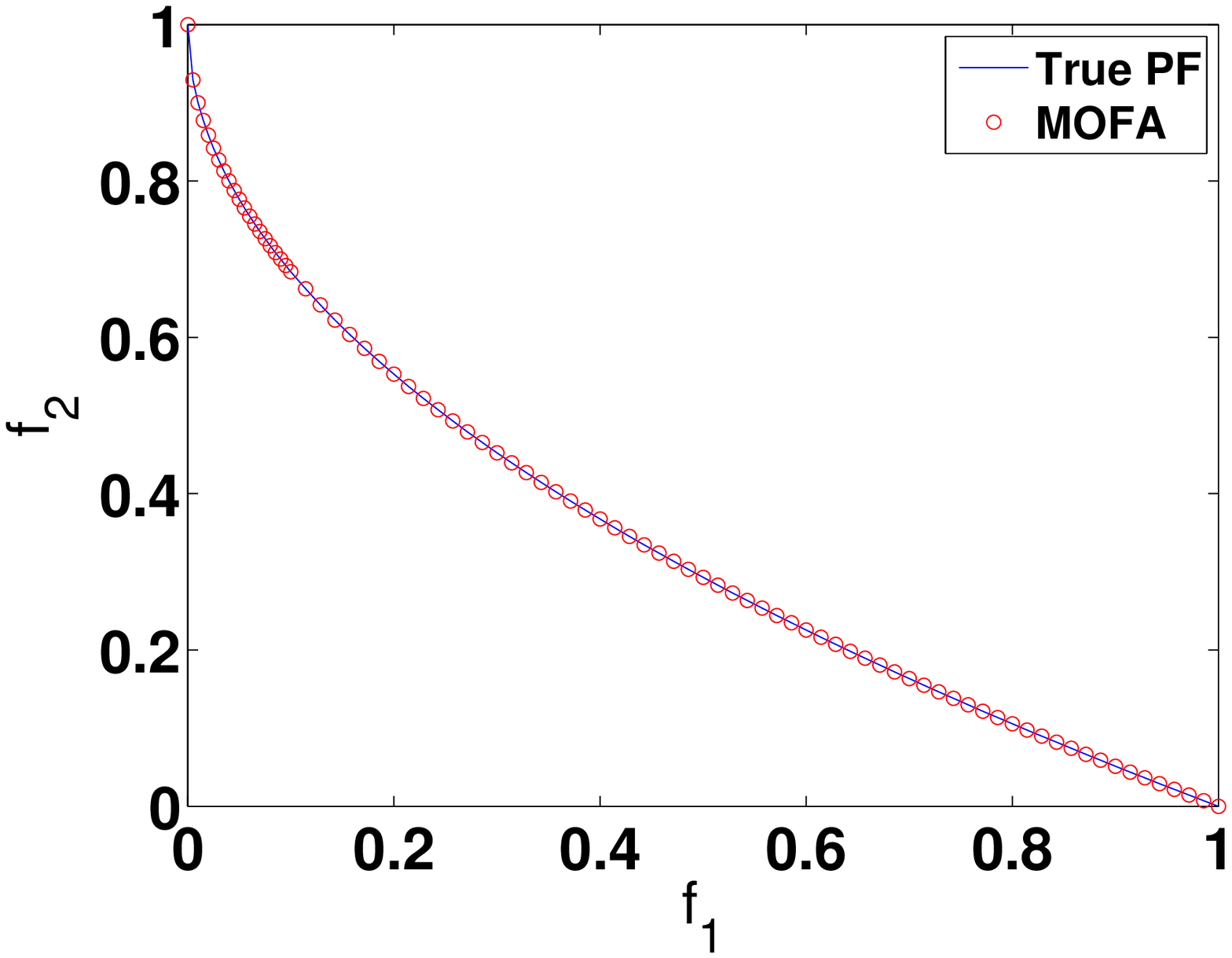} }
\caption{a) Pareto front of test function SCH, and b) Pareto front of test function ZDT2. \label{fig-300}}
\end{figure}

\begin{figure}
\centerline{\includegraphics[height=2.5in,width=3in]{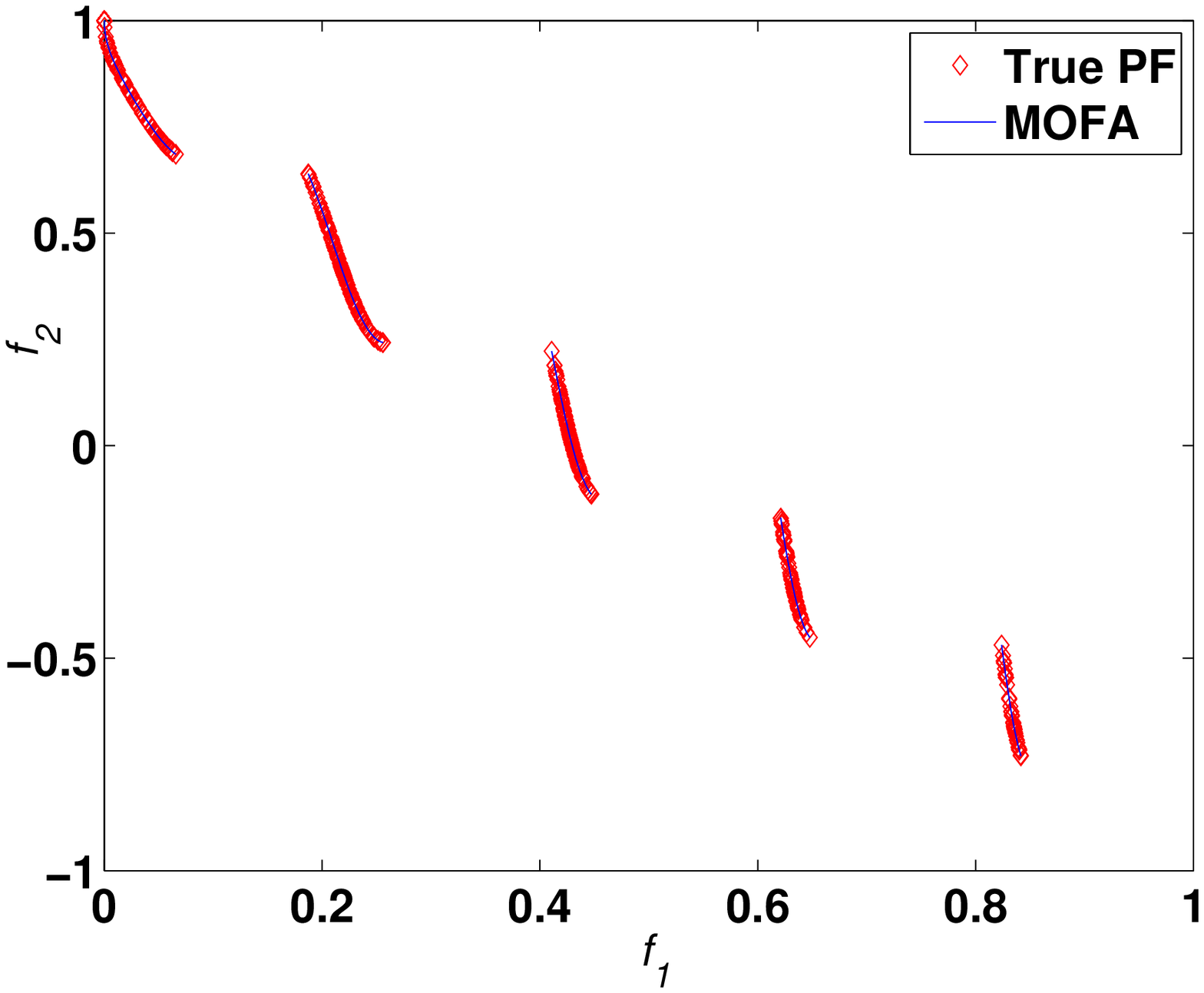}
\includegraphics[height=2.5in,width=3in]{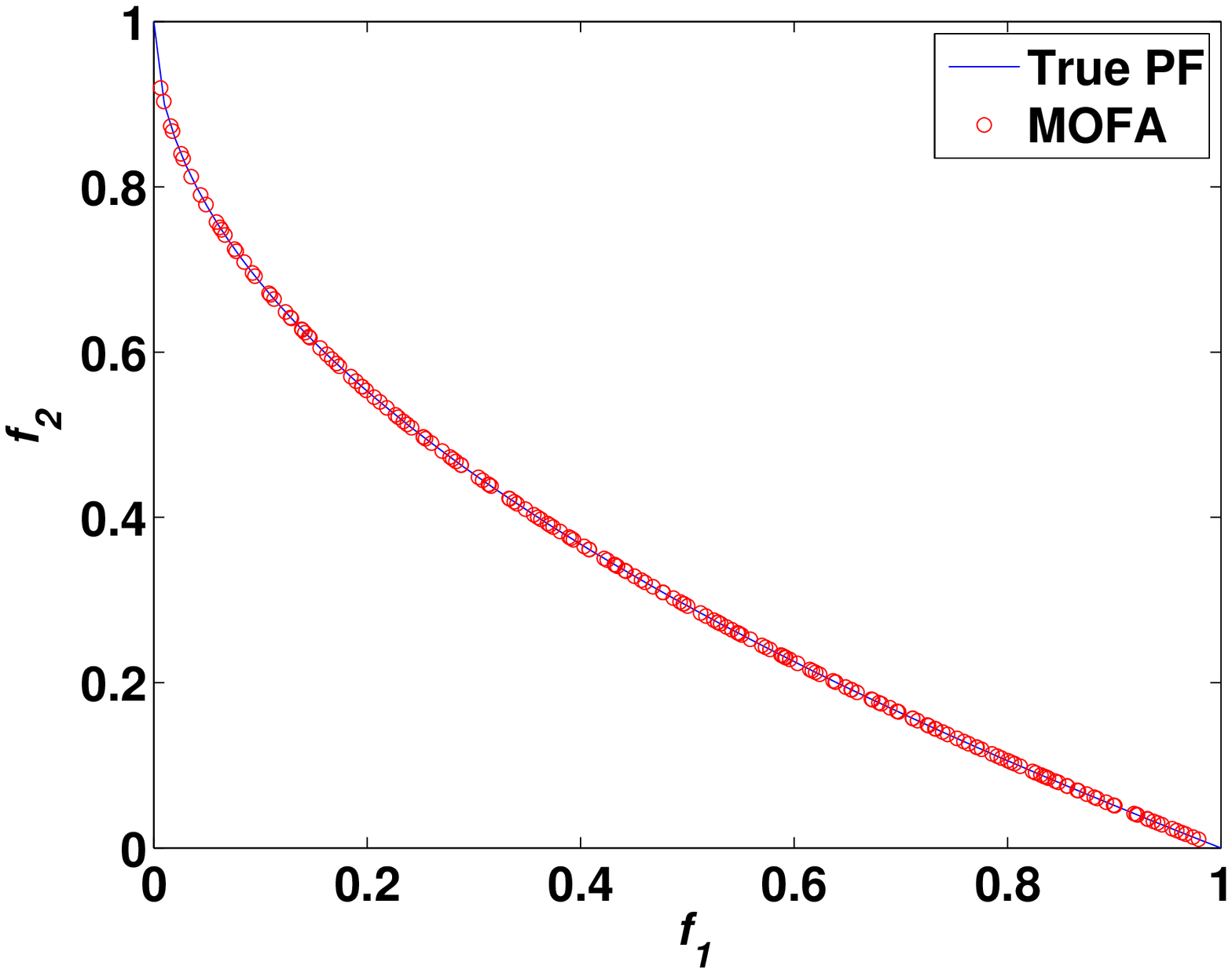} }
\caption{a) Pareto front of test function ZDT2, and b) Pareto front of test function LZ. \label{fig-400}}
\end{figure}

Fig. \ref{fig-200} shows the exponential-like decrease
of $D_g$ as the iterations proceed.
We can see clearly that our MOFA algorithm indeed converges almost exponentially.
The results for all the functions are summarized in Table 1, and the estimated
Pareto fronts and true fronts of other functions are shown in Fig. \ref{fig-300}
and Fig. \ref{fig-400}. In all these figures, the vertical axis is $f_2$ and
the horizontal axis is $f_1$.

\subsection{Comparison Study}

In order to compare the performance of the  proposed MOFA with other established
multiobjective algorithms, we have carefully selected a few algorithms with
available results from the literature. When the results have not been available,
we have implemented the algorithms using well-documented studies
and then generated new results using these algorithms. In particular, we have used other
methods for comparison, including vector evaluated genetic algorithm
(VEGA) \cite{Schaf}, NSGA-II \cite{Deb3}, multiobjective differential evolution (MODE) \cite{Xue,Babu},
differential evolution for multiobjective optimization (DEMO) \cite{Robic},
multiobjective bees algorithms (Bees) \cite{Pham},
and strength Pareto evolutionary algorithm (SPEA) \cite{Deb3,Mada}. The performance measures in terms of
generalized distance $D_g$ are summarized in Table 2 for all
the above major methods.

\begin{table}\caption{Summary of results.}
\begin{center} \begin{tabular}{|l|l|l|}
\hline
Functions & Errors (1000 iterations) & Errors (2500 iterations) \\
\hline
SCH & 5.5E-09  & 4.0E-22 \\
ZDT1 & 2.3E-6   & 5.4E-19  \\
ZDT2 & 8.9E-6   & 1.7E-14  \\
ZDT3 & 3.7E-5   & 2.5E-11  \\
LZ & 2.0E-6  & 7.7E-12 \\
\hline
\end{tabular} \end{center} \end{table}

It is clearly seen from Table 2 that the proposed MOFA obtained better
results for almost all five cases, though for ZDT2 function
our result is the same order
(still slightly better) as that by DEMO.

\begin{table}\caption{Comparison of $D_g$ for $n=50$ and $t=500$ iterations. }
\begin{center} \begin{tabular}{|l|l|l|l|l|l|l|}
\hline
Methods &  ZDT1 & ZDT2 & ZDT3 & SCH & LZ \\
\hline \hline
VEGA & 3.79E-02 & 2.37E-03  & 3.29E-01 & 6.98E-02 & 1.47E-03  \\
NSGA-II & 3.33E-02 & 7.24E-02  & 1.14E-01   & 5.73E-03 & 2.77E-02  \\
MODE & 5.80E-03 & 5.50E-03 &     2.15E-02 & 9.32E-04 & 3.19E-03  \\
DEMO & 1.08E-03 & 7.55E-04 &  1.18E-03 & 1.79E-04  & 1.40E-03 \\
Bees  & 2.40E-02 & 1.69E-02 & 1.91E-01  & 1.25E-02  & 1.88E-02 \\
SPEA & 1.78E-03 & 1.34E-03 & 4.75E-02 &5.17E-03  & 1.92E-03  \\ \hline
MOFA & 1.90E-04 & 1.52E-04 & 1.97E-04 & 4.55E-06 & 8.70E-04 \\
\hline
\end{tabular} \end{center} \end{table}

\section{Design Optimization}

Design optimization, especially design of structures, has many applications
in engineering and industry. As a result, there are many different benchmarks with
detailed studies in the literature \cite{Pham,Gando,Kim}. Some benchmarks have been solved
by various methods, while others do not have all available data for comparison. Thus, we have chosen
the welded beam design, and disc brake design among the well-known benchmarks
\cite{Kim,Rang,RAL}. In the rest of this paper,
we will solve these two design benchmarks using MOFA.

\subsection{Welded Beam Design}

Multiobjective design of a welded beam is a classical benchmark which has been solved by
many researchers \cite{Deb,Gong,RAL}. The problem has four design variables: the width $w$ and length $L$
of the welded area, the depth $d$ and thickness $h$ of the main beam. The objective is to minimize both
the overall fabrication cost and the end deflection $\delta$.

The detailed formulation can be found in  \cite{Deb,Gong,RAL}. Here we only rewrite the main problem as
\be \textrm{minimise } \; f_1(\x)=1.10471 w^2 L + 0.04811 d h (14.0+L),
\; \textrm{ minimize } \; f_2=\delta, \ee
subject to
\be
\begin{array}{lll}
g_1(\x)=w -h \le 0, \vspace{2pt} \\ \vspace{3pt}
g_2(\x) =\delta(\x) - 0.25 \le 0, \\ \vspace{3pt}
g_3(\x)=\tau(\x)-13,600 \le 0, \\ \vspace{3pt}
g_4(\x)=\sigma(\x)-30,000 \le 0, \\ \vspace{3pt}
g_5(\x)=0.10471 w^2 +0.04811  h d (14+L) -5.0 \le 0, \\ \vspace{3pt}
g_6(\x)=0.125 - w \le 0, \\ \vspace{3pt}
g_7(\x)=6000 - P(\x) \le 0,
\end{array}
\ee
where
\be \begin{array}{ll}
 \sigma(\x)=\frac{504,000}{h d^2},  & Q=6000 (14+\frac{L}{2}), \\ \\
 D=\frac{1}{2} \sqrt{L^2 + (w+d)^2}, & J=\sqrt{2} \; w L [ \frac{L^2}{6} + \frac{(w+d)^2}{2}], \\ \\
 \delta=\frac{65,856}{30,000 h d^3}, &  \beta=\frac{QD}{J}, \\ \\
 \alpha=\frac{6000}{\sqrt{2} w L}, & \tau(\x)=\sqrt{\alpha^2 + \frac{\alpha \beta L}{D}+\beta^2}, \\ \\
 P=0.61423 \times 10^6 \; \frac{d h^3}{6} (1-\frac{d \sqrt{30/48}}{28}). &
\end{array} \ee
The simple limits or bounds are $0.1 \le L, d \le 10$ and
$0.125 \le w, h \le 2.0$.

\begin{figure}
\centerline{\includegraphics[height=3in,width=4in]{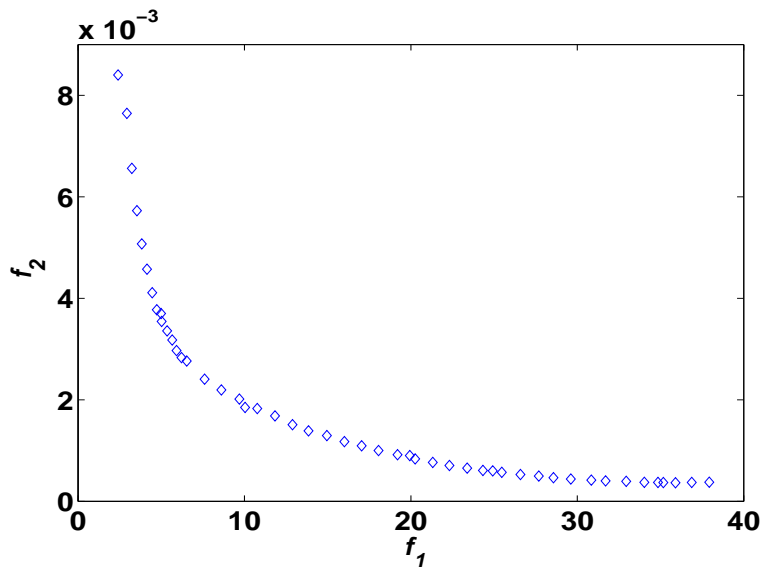} }
\caption{Pareto front for the bi-objective beam design where
the horizontal axis corresponds to cost and the vertical axis corresponds to deflection. \label{fig-440}}
\end{figure}

\begin{figure}
\centerline{\includegraphics[height=3in,width=4in]{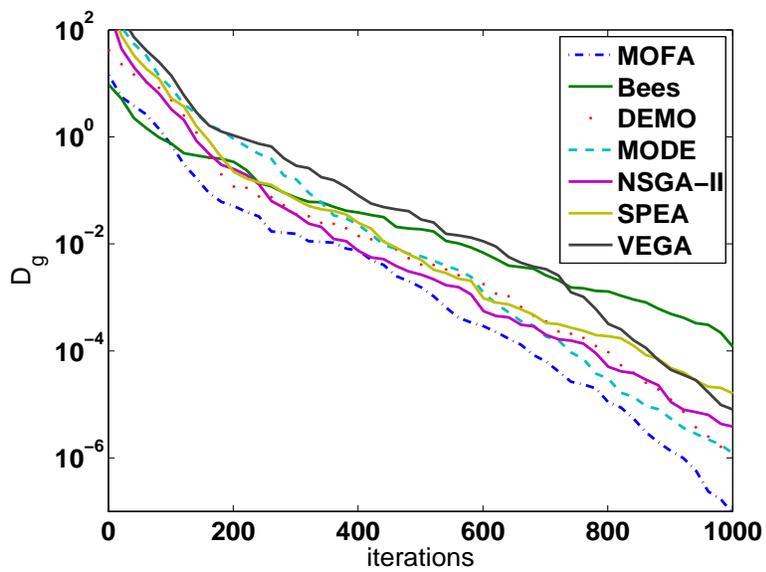} }
\caption{Convergence comparison for the beam design. \label{fig-800}}
\end{figure}

By using the MOFA, we have solved this design problem. The approximate Pareto front
generated by the 50 non-dominated solutions after 1000 iterations
are shown in Fig. \ref{fig-440}. This is
consistent with the results obtained by others \cite{Pham,RAL}. In addition, our results
are more smooth with fewer iterations, which shows the efficiency of the proposed MOFA.
A comparison of the results with those obtained by other methods is shown in Fig. \ref{fig-800}
where we can see MOFA converged faster.

\subsection{Design of a Disc Brake}

\begin{figure}[h]
\centerline{\includegraphics[height=3in,width=4in]{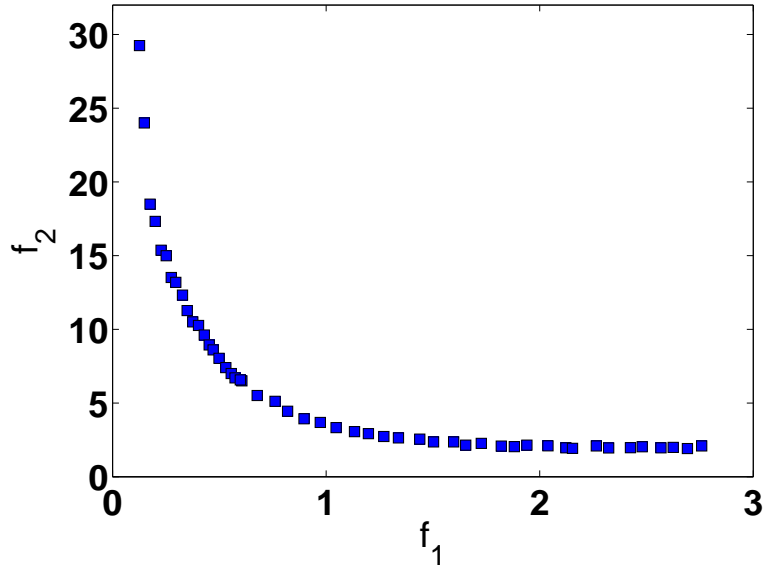} }
\caption{Pareto front for the disc brake design where $f_2$ is the vertical axis and $f_1$ is the horizontal axis. \label{fig-500}}
\end{figure}

\begin{figure}[h]
\centerline{\includegraphics[height=3in,width=4in]{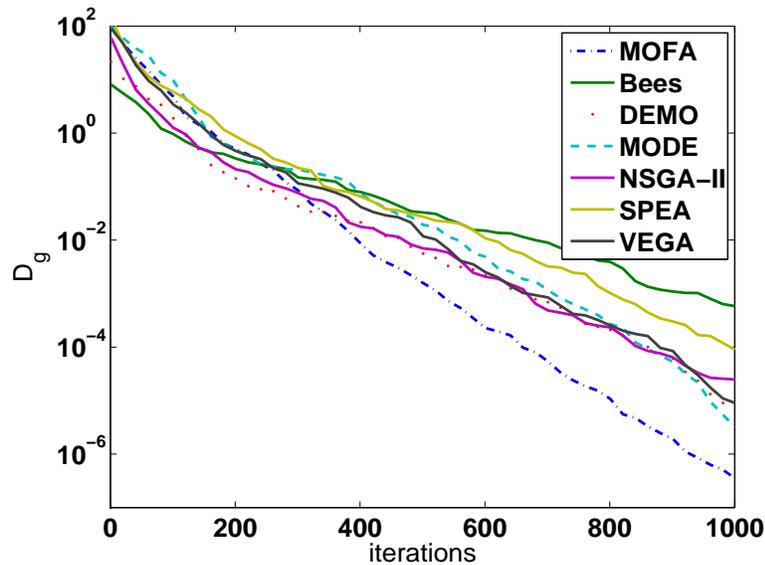} }
\caption{Convergence comparison for the disc brake design. \label{fig-900}}
\end{figure}

Design of a multiple disc brake is another benchmark for multiobjective optimization \cite{Gong,Osy,RAL}.
The objectives are to minimize the overall mass and the braking time by choosing optimal design
variables: the inner radius $r$, outer radius $R$ of the discs, the engaging force $F$ and the number
of the friction surfaces $s$. This is under the design constraints such as the torque, pressure, temperature,
and length of the brake. This bi-objective design problem can be written as:
\be \textrm{Minimize } \; f_1(\x) =4.9 \times 10^{-5} (R^2-r^2) (s-1),
\quad f_2(\x)=\frac{9.82 \times 10^6 (R^2-r^2)}{F s (R^3-r^3)}, \ee
subject to
\be
\begin{array}{lll}
g_1(\x) =  20-(R-r) \le 0, \\ \vspace{5pt}
g_2(\x) = 2.5 (s+1)-30 \le 0, \\ \vspace{5pt}
g_3(\x) = \frac{F}{3.14 (R^2-r^2)} -0.4 \le 0, \\ \vspace{5pt}
g_4(\x) = \frac{2.22 \times 10^{-3} F (R^3 -r^3)}{(R^2-r^2)^2} -1 \le 0, \\ \vspace{5pt}
g_5(\x) =900- \frac{0.0266 F s (R^3-r^3)}{(R^2-r^2)} \le 0.
\end{array}
\ee
The simple limits are \be 55 \le r \le 80, \; 75 \le R \le 110, \; 1000 \le F \le 3000, \; 2 \le s \le 20. \ee
The detail formulation of this problem and background descriptions can be found in \cite{Gando,Gong,Osy,RAL}.
The Pareto front of 50 solution points after 1000 iterations obtained by MOFA is shown in Fig. \ref{fig-500},
where we can see that the results are smooth and are the same or better than the
results obtained in \cite{RAL}. This can also be seen from the comparison of converge rates
shown in Fig. \ref{fig-900}.

In order to see how the proposed MOFA performs for the real-world design problems,
we also solved the same problems using other available multiobjective algorithms.
The comparison of the convergence rates is plotted in the logarithmic scales
in Fig. \ref{fig-800} and Fig. \ref{fig-900}. We can see from these two figures that
the convergence rates of MOFA are of the highest in an exponentially decreasing way in both cases.
This again suggests that MOFA provides better solutions in a more efficient way.

The simulations for these benchmarks and test functions suggest that MOFA is a very efficient
algorithm for multiobjective optimization. It can deal with highly nonlinear problems
with complex constraints and diverse Pareto optimal sets.

\section{Conclusions}

Multiobjective optimization problems are typically very difficult to solve.
We have successfully formulated a new algorithm for multiobjective
optimization, namely, multiobjective firefly algorithm (MOFA),
based on the recently developed firefly algorithm. The proposed MOFA has been tested against
 a subset of well-chosen test functions, and
then been applied to solve design optimization benchmarks in industrial engineering.

By comparing with other algorithms,
the present results suggest that MOFA is an efficient multiobjective optimizer.
Further studies can focus on parametric studies which
can be very useful to identify the optimal ranges
of parameters for various optimization problems.

In addition, convergence analysis of MOFA will also show insight in the working
mechanism of the algorithm, and this may also help to improve the proposed algorithm
or even design new algorithms.
For example,  hybridization with other algorithms may also prove to be fruitful.
Furthermore, formulation of a discrete MOFA will also be an important topic for further research.

\end{document}